\documentclass[12pt]{article}

\usepackage{amsfonts,amssymb,amsmath,amsthm,eucal}

\DeclareMathOperator{\Cov}{Cov}
\DeclareMathOperator{\Hur}{Hur}

\newcommand{\Z}{\mathbb{Z}}
\newcommand{\PP}{\mathbb{P}}
\newcommand{\fS}{\mathfrak{S}}
\newcommand{\la}{\lambda}
\newcommand{\be}{\beta}

\newcommand{\vac}{v_\emptyset}

\newcommand{\nr}[1]{:\!#1\!:}
\newcommand{\ul}{\underline}

\newtheorem*{Theorem}{Theorem}

\begin{document}

\title{Toda equations for Hurwitz numbers}
\author{Andrei Okounkov\thanks{
 Department of Mathematics, University of California at
Berkeley, Evans Hall \#3840, 
Berkeley, CA 94720-3840. E-mail: okounkov@math.berkeley.edu}
}
\date{} 
\maketitle

\begin{abstract} We consider ramified coverings of $\PP^1$
with arbitrary ramification type over $0,\infty\in\PP^1$ and
simple ramifications elsewhere and prove that the 
generating function for the numbers of
such coverings is a $\tau$-function
for the Toda lattice hierarchy of 
Ueno and Takasaki.
\end{abstract}

\section{Introduction}

\subsection{} 

In a recent paper \cite{P}, R.~Pandharipande showed that the conjectural Toda
equation for the Gromov-Witten potential $F_{\PP^1}$ of the Riemann sphere
$\PP^1$ implies a Toda equation for the generating function of the
Hurwitz numbers $H_{g,d}$. Here  $H_{g,d}$ is the number of 
degree $d$ connected 
coverings of $\PP^1$ by a curve of genus $g$ which have simple
 ramifications over some  $2g+2d-2$ fixed points of $\PP^1$. This Toda
equation results in a simple recurrence relation for the numbers $H_{g,d}$.

The purpose of this paper is to prove that a more general generating function,
where arbitrary ramifications are allowed over two given points of $\PP^1$
(for example, over $0,\infty\in\PP^1$)
satisfies an entire hierarchy of PDE's, namely the Toda lattice hierarchy of 
Ueno and Takasaki \cite{UT}. 
This proves and generalizes the result of Pandharipande.

We will call these more general numbers of ramified coverings the double
Hurwitz numbers. 

\subsection{}
Our approach is quite straightforward: we remark that the generating function
for double Hurwitz numbers is almost by definition a certain matrix
element in the infinite wedge space. It is a well known result of the Kyoto
school that such matrix elements are $\tau$-functions of integrable
hierarchies. 

These Toda equations should probably be not too surprising from 
the physical point of view. Physicists consider the enumeration 
of coverings of $\PP^1$ as a sort of string theory on $\PP^1$, 
see for example \cite{GT}, and relate it to a number of other models
where the occurrence of integrable hierarchies is common. 

\subsection{}
It is known, see \cite{ELSV} and also  \cite{GJV,GV}, that ordinary
Hurwitz numbers, that is, the numbers of coverings
where the ramification  over $\infty\in\PP^1$ is arbitrary and all other
ramifications are simple,
can be expressed as certain Hodge integrals over moduli spaces
of curves. It is an interesting problem to transfer the
Toda equation to the Hodge integrals side. In  particular, it would be 
nice to reproduce the KdV equations or Virasoro constraints,
 see for example \cite{K,Ge}, for just
the Gromov-Witten potential of a point, which corresponds
to Hodge integrals with no $\lambda$-classes.
 
It is also interesting and challenging to find an expression for the
double Hurwitz numbers 
in terms of the GW theory of $\PP^1$. 
The conjectural Toda equation for the GW potential of $\PP^1$,
see \cite{P}, may be an indication of such a connection.  

\subsection{}
The Toda lattice equations is an infinite number of PDE's
satisfied by the generating function for coverings, each
resulting is some recurrence 
relations for the double Hurwitz numbers. These
relations uniquely determine the double Hurwitz numbers. 

In particular, the lowest equation of the Toda hierarchy,
see Sections \ref{lt1} and \ref{lt2}, 
describes what will happen if we add an extra sheet 
to our covering leaving the nontrivial part of the 
monodromy the same. 

Many recurrent relations for Hurwitz numbers are known in the literature, see for
example \cite{GJV} and references therein, but as R.~Pandharipande writes 
in \cite{P}, the form of the
Toda recurrence relations is ``certainly the simplest recursion
for the Hurwitz numbers known to the author'' . It is absolutely
beyond the scope of the present note to give any sort
of survey on the vast literature on Hurwitz numbers and their applications.

\subsection{}
I want to thank Rahul Pandharipande for helpful correspondence. 

\section{Coverings and Toda equations} 

\subsection{} 

Let $S(d)$ denote
the symmetric group and fix some conjugacy classes 
$$
C_1,\dots,C_s\subset S(d)\,.
$$
Let $\Cov_d(C_1,\dots,C_s)$ denote the weighted  number of $d$-fold coverings of $\PP^1$ 
ramified over $s$ fixed points of $\PP^1$ 
with monodromies in the  conjugacy classes $C_1,\dots,C_s$.
The weight of any covering in $\Cov_d$ is the reciprocal of the
order of its group of automorphisms. This automorphisms group is the same as 
the centralizer of all monodromies in $S(d)$. Note that $\Cov_d$ counts all
coverings,
that is, including the disconnected ones.

For this weighted number there is the following formula, see for example \cite{D}, 
which goes back to Burnside 
\begin{equation}\label{BDW}
\Cov_d(C_1,\dots,C_s) = \sum_{|\la|=d} \left(\frac{\dim\la}{d!}\right)^2 \, \prod_{i=1}^s f_{C_i}(\la) \,,
\end{equation}
where the sum is over all partitions $\la$ of $d$, $\dim\la$ is the dimension of the 
corresponding representation of the symmetric group $S(d)$, and 
\begin{equation}\label{f}
f_{C_i}(\la) = |C_i|\, \frac{\chi^\la(C_i)}{\dim\la} \,.
\end{equation}
Here $|C_i|$ is the cardinality of the conjugacy class $C_i$ and $\chi^\la(C_i)$ is the value
of the irreducible character $\chi^\la$ on any element of $C_i$. The expression \eqref{f}
is in fact a polynomial in $\la$ in a certain natural sense, 
see for example \cite{OO}.

\subsection{} 

Recall that the Schur functions $s_\la$ are certain distinguished symmetric
functions  associated to 
any partition $\la$. In terms of the power-sum symmetric functions $p_1,p_2,\dots$ 
they have the following
expression, see \cite{M}, 
$$
s_\la(p_1,p_2,\dots) =\frac1{d!} \, \sum_{|\mu|=d} \chi^\la(C_\mu) \, |C_\mu|\, p_\mu \,,
$$
where the summation is over all partitions $\mu$ of $d=|\la|$, $C_\mu$ is the conjugacy class
corresponding to the partition $\mu$, and $p_\mu = \prod p_{\mu_i}$. We will treat 
the $p_k$'s   as independent variables. 

Now let $P=(p_1,p_2,\dots)$ and  
$P'=(p'_1,p'_2,\dots)$ be two sets of variables. We have the following equality: 
\begin{multline}\label{eq1}
\sum_{d=0}^\infty q^d \sum_{|\mu|=|\nu|=d}
p_\mu \, p'_\nu \, \Cov_d(C_\mu,C_\nu,C_1,\dots,C_s) 
=\\
\sum_\la q^{|\la|} \, s_\la(P)\, s_\la(P') \,  \prod_{i=1}^s f_{C_i}(\la)  \,,
\end{multline}
where $q$ is a new variable which keeps track of the degree $d$ of the covering. 
The variables $q$ is redundant since for every term its degree $d$ in $q$ is equal
to both its degree in $P$ and in $P'$, where  
$$
\deg p_k = \deg p'_k = k \,,
$$
but it is convenient to keep it. 

\subsection{} 

Let $C_{(2)}$ be the conjugacy class of a transposition.   Introduce the following 
generating function
$$
\tau(P,P',\be,q)=\sum_{d,b,\mu,\nu}  q^d \, \be^b \, p_\mu \, p'_\nu \, 
\Cov_d(C_\mu,C_\nu,\underbrace{C_{(2)},\dots,C_{(2)}}_{\textup{$b$ times}})
\big/b! \,.
$$
This is the generating function for degree $d$ possibly disconnected 
coverings of $\PP^1$ such that the 
ramification over $0,\infty\in\PP^1$ 
 can be fixed arbitrarily and the other $b$ ramifications
are simple. It is a general property of the exponential generating functions,
see e.g.\ Chapter 3 in \cite{GJ},  that
$$
\log \tau(P,P',\be,q) =  H(P,P',\be,q)\,,
$$
where $H(P,P',\be,q)$ is the generating function for connected coverings
$$
H(P,P',\be,q)= \sum_{d,b,\mu,\nu}  q^d \, \be^b \, p_\mu \, p'_\nu \, 
\Hur_{d,b}(\mu,\nu)\big/b!\,,
$$ 
Here $\Hur_{d,b}(\mu,\nu)$ is the weighted number of \emph{connected} degree $d$ 
coverings of $\PP^1$ with monodromy around $0,\infty\in\PP^1$ being $\mu$ and $\nu$,
respectively, and $b$ 
additional simple ramifications. The genus of such a covering is
$$
g = (b+2-\ell(\mu)-\ell(\nu))/2\,,
$$
where $\ell(\mu)$ is the number of parts of $\mu$. The weight of each covering in
$\Hur_{d,b}(\mu,\nu)$ is the reciprocal of the order of its automorphism
group. We call the numbers $\Hur_{d,b}(\mu,\nu)$ the \emph{double Hurwitz 
numbers}.

We also point out that if when one of ramifications $\mu$ and $\nu$
is trivial and the covering is connected of degree $d>2$ 
then the automorphism group is
trivial and so $H=\log\tau$ counts any  such covering with weight $1$. 
 
\subsection{}

Denote by $f_2$ the polynomial \eqref{f} corresponding to the class $C_{(2)}$ of
a transposition.  From \eqref{eq1} we have
\begin{equation}\label{eq2}
\tau(P,P',\be,q)=\sum_\la q^{|\la|} \,  e^{\be\,f_2(\la)} \,  s_\la(P)\, s_\la(P')   \,. 
\end{equation}
Explicitly, the polynomial $f_2$ is given by  
\begin{equation}\label{f21}
f_2=\frac12 \sum_{i} \left[\big(\la_i-i+\tfrac12\big)^2 - \big(-i+\tfrac12\big)^2\right] \,.
\end{equation}
Remark that the upper limit of summation can be
taken to be $+\infty$ in the above sum since $\la_i=0$ for all sufficiently
large $i$. Using some standard combinatorics, see for example Section 5 of \cite{BO}, one can rewrite 
the formula \eqref{f21} as follows
\begin{equation}\label{f22}
f_2=\sum_{k\in \fS(\la)_+} \frac{k^2}2 - \sum_{k\in \fS(\la)_-} \frac{k^2}2 
\end{equation}
where 
$$
\fS(\la)=\{\la_i-i+\tfrac12\} \subset \Z+\tfrac12 
$$
and, by definition, for any $S\subset\Z+\frac12$ we set 
$$
S_+ = S \setminus \left(\Z_{\le 0} - \tfrac12 \right) \,, \quad
S_- = \left(\Z_{\le 0} - \tfrac12 \right) \setminus S \,. 
$$

\subsection{} 

The sum \eqref{eq2} admits a straightforward interpretation as a certain matrix
element in the infinite wedge space. It will be convenient for us to use
the notation of \cite{iw}, see especially the Appendix to \cite{iw} for a summary
of the infinite wedge basics. A good general reference on the infinite wedge space
is Chapter 14 of \cite{K}. 

Note that the power-sum  variables 
 $P$ and $P'$ are related to the variables $t$ and $t'$ used in
\cite{iw} by 
$$
t_k = \frac{p_k}{k} \,, \quad t'_k = \frac{p'_k}{k}\,, \quad k=1,2,\dots\,, 
$$
We will use the following abbreviations
$$
\Gamma_+ = \Gamma_+\left(p_1,\frac{p_2}{2},\frac{p_3}{3},\dots\right)\,, \quad 
\Gamma_- = \Gamma_-\left(p'_1,\frac{p'_2}{2},\frac{p'_3}{3},\dots\right)\,. 
$$
The formula \eqref{eq2} becomes 
\begin{equation}\label{eq3}
\tau(P,P',\be,q)=\left(\Gamma_+\, q^{H}e^{\be F_2} \, \Gamma_-\,  \vac, \vac\right) \,,
\end{equation}
where $\vac$ is the vacuum vector of the infinite wedge space, 
$H$ is the energy operator, and $F_2$ is the following
operator
$$
F_2= \sum_{k\in\Z+\frac12} \frac{k^2}2 \, \nr{\psi_k\, \psi_k^*} \,,
$$
where the colons denote the normal ordering. It is clear from \eqref{f22} that
the operator $F_2$ acts as follows
$$
F_2 \, v_\la = f_2(\la) \, v_\la \,,
$$
where 
$$
v_\la = \ul{s_1} \wedge \ul{s_2} \wedge \ul{s_3} \wedge \dots\,,
\quad \{s_1,s_2,\dots\}=\fS(\la)\,,
$$ 
is the vector corresponding to the partition $\la$. 

\subsection{}

Any operator of the form $\xi=\psi_i \psi^*_j$ satisfies the relation
\begin{equation*}
\left[\xi \otimes 1 + 1 \otimes \xi , 
\Omega\right] = 0 \,, \quad \Omega=\sum \psi_k \otimes \psi^*_k \,,
\end{equation*}
Since the operators  $q^{H}$ and $e^{\be F_2}$ are exponentials of sums
of such operators $\xi$, we have 
\begin{equation*}
\left[q^{H}e^{\be F_2} \otimes q^{H}e^{\be F_2}, 
\Omega\right] = 0 \,. 
\end{equation*}
It follows that the sequence 
\begin{equation}\label{taun}
\tau_n = \left(\Gamma_+\, q^{H}e^{\be F_2} \, \Gamma_-\,  v_n, v_n\right)\,, \quad n\in\Z\,,
\end{equation}
where 
\begin{equation*}
v_n = \ul{n-\tfrac12} \wedge \ul{n-\tfrac32} \wedge \ul{n-\tfrac52} \wedge \dots\,,
\end{equation*}
is the vacuum vector in the charge $n$ subspace,
is a sequence of $\tau$-functions for the Toda lattice hierarchy. Our
old  $\tau(P,P',\be,q)$ is the $\tau_0$ term of this sequence. 

\subsection{}

We will now show that all terms of the sequence \eqref{taun} can be
expressed in terms of $\tau(P,P',\be,q)$. We have  $v_n = R^n\, \vac$, where
$R$ is the translation operator, and also $[R,\Gamma_{\pm}]=0$.  We compute
\begin{align*}
R^{-n} F_2 R^n  & = \sum_{k>0} \frac{k^2}2 \, \psi_{k-n} \, \psi^*_{k-n} -
\sum_{k<0} \frac{k^2}2 \, \psi^*_{k-n} \, \psi_{k-n}  \\
& = \sum_{k\in\Z+\frac12} \frac{(k+n)^2}2 \, \nr{\psi_{k} \, \psi^*_{k}}  + 
\sum_{k=1/2}^{n-1/2} \frac{k^2}{2} \\
& = F_2 + n \,H + \frac{n^2}2\, C + \frac{n(4n^2-1)}{24} \,,
\end{align*}
where $H$ and $C$ are the energy and charge operators. Similarly, 
$$
R^{-n} \, H \, R^n = H+nC+\frac{n^2}2 \,.
$$
Since the charge operator $C$ commutes with $\Gamma_\pm$ and $C\vac=0$, we have
\begin{align*}
\tau_n &= \left(\Gamma_+\, q^{H}e^{\be F_2} \, \Gamma_-\,  v_n, v_n\right) \\
& =  \left(\Gamma_+\, R^{-n} q^{H}e^{\be F_2}  R^n \, \Gamma_-\,  \vac, \vac\right)\\
& =  q^{n^2/2} e^{n(4n^2-1)\be/24} \, 
\left(\Gamma_+\,  (e^{n\be} q)^{H}\, e^{\be F_2}   \, \Gamma_-\,  \vac, \vac\right)\\
& =  q^{n^2/2} e^{n(4n^2-1)\be/24} \, \tau(P,P',\be,e^{n\be} q) \,.
\end{align*}
Substituting this expression into the Hirota bilinear equations for
the Toda lattice hierarchy,
we obtain the following

\begin{Theorem} The function $\tau(P,P',\be, q)$ satisfies the
following Hirota-type bilinear equations: 
\begin{equation}\label{Hirota}
q^{m+1} e^{m(m+1)\be/2} [z^{-1-m}]\,  e^{-2\sum n \, s_n/z^n} 
\,\diamondsuit\, \heartsuit  = [z^{m+1}]\, e^{2\sum n \, s'_n/z^n}\,  
\clubsuit\, \spadesuit \,,
\end{equation}
where $m\in\Z$,  $z$ is a formal variable, $[z^{m+1}]$
denotes the coefficient of $z^{m+1}$, $S=(s_1,s_2,\dots)$ and $S'=(s'_1,\dots)$
are two arbitrary sequences, 
\begin{alignat*}{2}
\diamondsuit&=\tau(P+S,P'+S'+\vec z,\be,e^{(m+1)\be}q), \,\,\, & 
\clubsuit&=\tau(P+S-\vec z,P'+S',\be,e^{m\be}q) \,,  \\
\heartsuit&=\tau(P-S,P'-S'-\vec z,\be,e^{-\be}q)\,,  \quad &
\spadesuit&=\tau(P-S+\vec z,P'-S',\be,q)\,,
\end{alignat*}
and $\vec z = (z,z^2,z^3,\dots)$ \,. 
\end{Theorem}

Remark that, as pointed out earlier, 
\begin{align*} 
\tau(P,P',\be, q)&=\tau(q\, p_1, q^2\, p_2, q^3 \, p_3,\dots, P',\be,1)\\
&=\tau(P,q \,p'_1, q^2 \,p'_2, q^3 \, p'_3,\dots,\be,1)\,.
\end{align*}

\subsection{}\label{lt1}

In particular, setting $m=0$ and taking the coefficient of $s_1$ 
gives 
$$
\tau \, \frac{\partial^2\tau}{\partial p_1  \partial p'_1} -
\frac{\partial \tau}{\partial p_1} \, \frac{\partial \tau}{\partial p'_1}   =
q\, \tau(\dots, e^\be q) \, \tau(\dots, e^{-\be}q)\,.
$$  
This is equivalent to
\begin{equation}\label{toda}
\frac{\partial^2}{\partial p_1  \partial p'_1}
\log \tau = q \,  \frac{\tau(\dots, e^\be q) \, \tau(\dots, e^{-\be}q)}{\tau^2} \,.
\end{equation}
Recall that 
$$
\log\tau=H
$$ 
is the generating function for the
connected coverings.

The combinatorial meaning of the operator 
$\dfrac{\partial^2}{\partial p_1  \partial p'_1}$
is that it peels off a sheet of the covering which is
unramified over the special points $0,\infty\in\PP^1$.
The equation \eqref{toda} gives thus a recursive procedure
for the computation of the double Hurwitz numbers with a 
fixed nontrivial part of the  monodromy around  $0,\infty$.  

\subsection{}\label{lt2} 

In the particular case when there is no ramification over $0,\infty$,
that is, when,
$$
p_2=p_3=\dots=p'_2=p'_3=\dots=0
$$
the function $\tau$ depends on the variables $p_1$, $p'_1$, and $q$ 
only via their product $q p_1 p'_1$. If we set 
$$
u = \log q \, p_1 \, p'_1
$$
then the equation \eqref{toda} becomes
$$
e^{-u} \, \frac{\partial^2}{\partial u^2} \, \log \tau =
\frac{\tau(u+\be,\be)\, \tau(u-\be,\be)}{\tau^2} \,,
$$
which is equivalent to the Toda equation derived in \cite{P}.


\begin{thebibliography}{99} 

\bibitem{BO}
S.~Bloch and A.~Okounkov,
\emph{The character of the infinite wedge representation},
Adv.\ Math.\ \textbf{149} (2000), 1--60, 
alg-geom/9712009.  

\bibitem{D}
R.~Dijkgraaf,
\emph{Mirror symmetry and elliptic curves},
The Moduli Space of Curves, R.~Dijkgraaf,
C.~Faber, G.~van~der~Geer (editors), 
Progress in Mathematics, \textbf{129},
Birkh\"auser, 1995.

\bibitem{ELSV}
T.~Ekedahl, S.~Lando, M.~Shapiro, and A.~Vainshtein,
\emph{On Hurwitz numbers and Hodge integrals},
math.AG/9902104.  

\bibitem{GV} 
T.~Graber and R.~Vakil,
\emph{Hodge integrals and Hurwitz numbers via
virtual localization},
math.AG/0003028 

\bibitem{Ge}
E.~Getzler,
\emph{The Virasoro conjecture for Gromov-Witten invariants},
Algebraic geometry: Hirzebruch 70 (Warsaw, 1998), 147--176, Contemp.\ Math., 241,
Amer.\ Math.\ Soc., 1999. 

\bibitem{GJ}
I.~P.~Goulden and D.~M.~Jackson,
\emph{Combinatorial enumeration},
John Wiley \& Sons, 1983.

\bibitem{GJV} 
I.~P.~Goulden, D.~M.~Jackson, and R.~Vakil,
\emph{The Gromov-Witten potential of a point, Hurwitz
numbers, and Hodge integrals},
math.AG/9910004. 

\bibitem{GT} 
D.~J.~Gross and W.~Taylor,
\emph{Two-dimensional QCD and strings},
hep-th/9311072.


\bibitem{K}
V.~Kac,
\emph{Infinite dimensional Lie algebras},
Cambridge University Press.

\bibitem{Ko}
M.~Kontsevich,
\emph{Intersection theory on the moduli space of curves 
and the matrix Airy function}, Comm.\ Math.\ Phys.\ \textbf{147}
(1992), no.~1, 1--23. 

\bibitem{M}
I.~G.~Macdonald,
\emph{Symmetric functions and Hall
polynomials},
Clarendon Press, 1995. 

\bibitem{iw}
A.~Okounkov,
\emph{Infinite wedge and random partitions},
math.RT/9907127.


\bibitem{OO}
A.~Okounkov and G.~Olshanski,
\emph{Shifted Schur functions}
Algebra i Analiz \textbf{9}, 1997, no.~2, 73--146; translation in St.\
Petersburg Math.\ J.\ \textbf{9}, 1998, no.~2, 239--300. 


\bibitem{P}
R.~Pandharipande,
\emph{The Toda equations and the Gromov-Witten theory of the Riemann
sphere},
math.AG/9912166. 

\bibitem{UT}
K.~Ueno and K.~Takasaki,
\emph{Toda lattice hierarchy},
Adv.\ Studies in Pure Math.\ \textbf{4}, Group Representations
and Systems of Differential Equations, 1--95, 1984. 



\end{thebibliography}
\end{document}